\documentclass[twoside,a4paper,reqno,12pt]{amsart}

\usepackage{cite}
\usepackage{mathrsfs}
\usepackage{amsmath}
\usepackage{latexsym}
\usepackage{amssymb,bm}
\usepackage{amsfonts}
\usepackage{longtable}
\usepackage{multirow, hyperref}
\usepackage{braket}
\usepackage{graphicx}
\usepackage{cleveref}
\usepackage{amsmath,enumerate}
\usepackage{tikz}
\usetikzlibrary{calc}
\usetikzlibrary{decorations.markings, arrows.meta}

\usepackage{graphics}

\usepackage{color}

\def\Frob{{\rm Fro}}
\def\Irr{{\rm Irr}}
\def\Syl{{\rm Syl}}
\def\Sol{{\rm Sol}}
\def\Fit{{\rm Fit}}
\def\cod{{\rm cod}}

\newtheorem{thm}{Theorem}[section]

\newtheorem{proposition}[thm]{Proposition}

\newtheorem{defi}[thm]{Definition}

\newtheorem{lemma}[thm]{Lemma}

\newtheorem{corollary}[thm]{Corollary}

\theoremstyle{definition}

\def\leq{\leqslant}

\def\geq{\geqslant}

\begin{document}

 	\title{On the Codegree graphs of finite groups}

	\author{ Jiyong Chen}
	\author{Ni Du}
	\author{Leyi Li}

	\address{J.  Chen, School of Mathematical Sciences\\
		Xiamen University \\
		Xiamen 361005\\
		P. R. China}
	\email{chenjy1988@xmu.edu.cn}
	\address{N. Du, School of Mathematical Sciences\\
		Xiamen University \\
		Xiamen 361005\\
		P. R. China}
	\email{duni@xmu.edu.cn}
	\address{L. Li, School of Mathematical Sciences\\
		Xiamen University \\
		Xiamen 361005\\
		P. R. China}
	\email{lileyi@stu.xmu.edu.cn}

	\date\today

\begin{abstract}
The codegree of an irreducible character $\chi$ of a finite group $G$ is defined as $|G:\ker\chi|/\chi(1)$. The codegree graph $\Gamma(G)$ of a finite group $G$ is the graph whose vertices are the prime divisors of $|G|$, where two distinct primes $p$ and $q$ are adjacent if and only if $pq$ divides the codegree of some irreducible character of $G$. In this paper, we prove that a graph can occur as a codegree graph $\Gamma(G)$ of some finite group $G$ if and only if its complement is triangle-free and $3$-colorable. This generalizes the known characterization for codegree graphs from solvable groups to all finite groups.
As an application, we give a full classification of all groups for which $\Gamma(G)$ is a $5$-cycle. We also investigate conditions under which the codegree graph coincides with or differs from the prime graph for solvable groups.

\end{abstract}
		\maketitle

 \section{Introduction}
 Throughout this paper, all groups are assumed to be finite, all graphs refer to finite simple graphs, and all characters are considered as complex characters.
 
  The \emph{codegree} of an irreducible character $\chi$ of a group $G$ is defined in \cite{qian2002} as $$\cod(\chi)=\frac{|G:ker\chi|}{\chi(1)}.$$ The set of  codegrees of a group $G$ is denoted by $\cod(G)$. Due to the natural connection with character degrees, numerous results have been obtained concerning character codegrees, see \cite{alizadeh2019} \cite{du2016} \cite{qian_overview} \cite{yang2017}  .
 
 The \emph{Gruenberg-Kegel graph} of a set of integers $N$ is defined as follows. Its vertices are the prime numbers that divide some element of $N$, and two distinct vertices $p$ and $q$ are adjacent if and only if there exists an integer $n \in N$ such that $pq$ divides $n$. 
 The \emph{codegree graph} of a group $G$, denoted by $\Gamma(G)$, is the Gruenberg-Kegel graph of the set $\cod(G)$. The concept of codegree graph was introduced by Qian in \cite{qian2007}. 
The codegree graph can reveal several properties of the group $G$. 
 For example, if $G$ is unsolvable, then its codegree graph $\Gamma(G)$ contains a triangle\cite{alizadeh2021}. Furthermore, if the Fitting subgroup of $G$ is trivial, then the codegree graph $\Gamma(G)$ must be a complete graph, see \cite{alizadeh2020}.
 
 A natural question arises: which graphs can occur as codegree graphs? To address this question, we employ the prime graph as a tool.
 The \emph{prime graph} of a group $G$, denoted by $\Gamma_e(G)$, is defined as the Gruenberg–Kegel graph of the set of element orders of $G$.
In \cite{qian2007}, Qian has proven that the prime graph $\Gamma_e(G)$ of any group $G$ is a subgraph of its codegree graph $\Gamma(G)$. Therefore, the properties of prime graphs are closely related to those of codegree graphs, and we can use this relationship to investigate codegree graphs.

A notable result on the structure of prime graph was established by Gruber, Keller and Lewis. They provided a characterization of the prime graphs for solvable groups in \cite{gruber2015}. Then Liu and Yang \cite{liu2025} obtained a parallel result for the codegree graphs of solvable groups. 
In this paper, we extend Liu and Yang's work from solvable groups to all finite groups. We determine precisely which graph can occur as the codegree graph of finite groups in the following theorem. 

 \begin{thm}\label{thm:codegree_graph}
A graph $\Gamma$ can occur as the codegree graph of some finite group $G$ if and only if its complement $\overline{\Gamma}$ is both triangle-free and $3$-colorable. 
\end{thm}
What surprises us the most is that a graph can occur as the codegree graph of some group if and only if it can occur as the codegree graph of some solvable group.

\Cref{thm:codegree_graph} shows that the complement $\overline{\Gamma}$ of the codegree graph of a group $G$ is always triangle-free. In fact, if the codegree graph itself is triangle-free, it can reveal many properties of the group $G$.
 By \cite{alizadeh2021}, groups with triangle-free codegree graphs are solvable. In \cite[Question 3.1]{qian_overview}, Qian asked for a characterization of those groups. 
 Classical Ramsey theory states that if a graph and its complement are both triangle-free, then the graph has at most $5$ vertices.
 The only $5$-vertex graph such that both the graph and its complement are triangle-free is the $5$-cycle.
In this paper, 
we give a complete classification of the finite groups whose codegree graph is a $5$-cycle.

\begin{thm}\label{thm:$5$-cycle}
  Let $G$ be a finite group. Suppose that $\pi(G)=\{a,b,c,d,e\}$. Let $A,B,C,D,E$ be Sylow $a$-, $b$-, $c$-, $d$-, $e$-subgroups of $G$ respectively, such that $A,B,C$, $D,E$ form a Sylow system of $G$.  Then, up to a permutation of the set $\pi(G)$, the codegree graph $\Gamma(G)$ is a $5$-cycle if and only if one of the following conditions holds. 

\begin{itemize}
        \item [(1)] $G=\Fit(G) \rtimes (B\rtimes 
        AD)$, where $\Fit(G)=C \times E$, $ABC=\Frob_2(A,B,C)$, $DCE=\Frob(D,CE)$, $AE= \Frob(A,E)$.
        \item [(2)] $G=\Fit(G) \rtimes (EB\rtimes \overline{A}D)$, where $\Fit(G)=C\times O_a(G)$, $\overline{A}\cong A/O_a(G)>1$, $ABC=\Frob_2 (A,B,C)$, $DCE=\Frob(D,CE)$, $AEB=\Frob_2(\overline{A},BE,O_a(G))$.
   \end{itemize}
   In particular, for a solvable group $G$, the codegree graph $\Gamma(G)$ is a $5$-cycle if and only if its prime graph $\Gamma_e(G)$ is a $5$-cycle.  
\end{thm}

Furthermore, we explore conditions under which the codegree graph and the prime graph of a solvable group are identical or distinct.

 \section{Preliminaries}

Frobenius groups and 2-Frobenius groups are important tools in this paper for studying the codegree graph. We begin by recalling their definitions.
\begin{defi}
A group $G$ is called a \emph{Frobenius group} if it satisfies $G = N \rtimes A$, where $A$ acts on $N$ such that $C_N(a) = 1$ for every non-identity element $a \in A$. Here, $N$ is the Frobenius kernel, $A$ is the Frobenius complement, and we call that the action of $A$ on $N$ is fixed-point-free. We denote such a Frobenius group $G$ by $\Frob(N,A)$.
\end{defi}
\begin{defi}
A group $G$ is called a \emph{2-Frobenius group} if there exists normal subgroups $M$ and $N$ of $G$ such that $M$ is a Frobenius group with kernel $N$, and $G/N$ is a Frobenius group with kernel $M/N$. In particular, if $N$ is a $p$-group, $M/N$ is a $q$-group, and $G/M$ is an $r$-group, then $G$ is denoted as a 2-Frobenius group of type $(p, q, r)$.
We denote such a 2-Frobenius group $G$ by $\Frob_2(M,N,A)$.
\end{defi}

The relationship between the codegree graph and Frobenius groups is connected through the Hall subgroups corresponding to pairs of non-adjacent primes in $\Gamma(G)$.
 \begin{lemma}\label{lem:disconnection_to_hall_subgroup}
\cite{qian2007} Let $pq \in \overline{\Gamma}(G)$, then $G$ admit a Hall $\{p,q\}$-subgroup $PQ$. And either $PQ$ is a Frobenius group or a 2-Frobenius group.
\end{lemma}
Moreover, the Frobenius kernel of the Hall subgroup from Lemma \ref{lem:disconnection_to_hall_subgroup} lies in the solvable radical of $G$.
\begin{lemma}\label{lem:Q_leq_Sol(G)}
    \cite{qian2025} Let $pq \in \overline{\Gamma}(G)$, if the Hall $\{p,q\}$-subgroup $PQ$, where $P$ and $Q$ are respectively Sylow $p$- and $q$- subgroup, satisfying one of the following properties:
    \item [(1)] $PQ$ is a Frobenius group with kernel $Q$ and complement $P$;
    \item [(2)] $PQ$ is a $(p,q,p)$-type 2-Frobenius group;
    then $Q \leq \Sol(G)$, where $\Sol(G)$ denotes the largest solvable normal subgroup of $G$. 
\end{lemma}
The special case when the Frobenius kernel is normal enables us to investigate the prime divisors of the codegrees in $G$.
\begin{lemma}\label{lem:normal_sylow634}
    Let $P$ be a Sylow $p$-subgroup of $G$, and suppose that $P$ acts fixed-point-freely on a normal subgroup $N$ of $G$. For any irreducible character $\chi \in \Irr(G)$, if $N\nleq ker \chi$, it follows that $p \nmid \cod(\chi)$.
\end{lemma}

\begin{proof}
    Let $\chi \in \Irr(G)$ be an irreducible character such that $N \nleq ker\chi$. 
    Then the restriction $\chi_N$ is not a multiple of the principal character of $N$, and thus it has at least one nonprincipal irreducible constituent $\lambda\in \Irr(N)$. By \cite[Theorem 6.34]{isaac_character}, since $P$ acts fixed-point-freely on $N$, it follows that $P \cap I_G(\lambda)=1$. Since $I_G(\lambda)\trianglelefteq G$, and $P$ is a Sylow $p$-subgroup of $G$, $p\nmid |I_G(\lambda)|$. By \cite[Theorem 6.11]{isaac_character}, there exists an irreducible character $\psi \in \Irr(I_G(\lambda))$ such that $\psi^G=\chi$. Therefore, 
    $$\cod(\chi)=\frac{|G:ker\chi|}{\chi (1)}=\frac{|G:ker\chi|}{\psi (1)|G:I_G(\lambda)|}=\frac{|I_G(\lambda)|}{\psi (1)|ker\chi|}$$
    It follows that $p\nmid \cod(\chi)$. 
\end{proof}
\begin{corollary}
    Let $P$ be a normal Sylow $p$-subgroup of $G$. Then the vertex $p$ has the same adjacency relation in the codegree graph as in the prime graph.
\end{corollary}

\section{Characterization of Codegree Graphs}

In this section, we characterize the codegree graphs of every finite groups. We extend the work of Liu and Yang \cite{liu2025}, who characterized these graphs for solvable groups, to the general case. Specifically, a graph is realizable as the codegree graph of some finite group $G$ if and only if its complement $\overline{\Gamma}$ is both triangle-free and $3$-colorable. 
The proof adapts the concept of the Frobenius digraph from \cite{gruber2015} to analyze adjacency in the codegree graph.

Extending results from solvable to arbitrary finite groups is often difficult, because Hall subgroups may not exist. However, we do not need to analyze simple groups case by case when we consider codegree graphs of non-solvable groups. This is because for any non-adjacent pair of vertices $p,q$ in $\Gamma(G)$, by Lemma~\ref{lem:Q_leq_Sol(G)}, the required Hall $\{p,q\}$-subgroup already exists in $\Sol(G)$. That is a key observation that simplifies our proof. 

We can define an orientation on the edges of the complement codegree graph $\overline{\Gamma}(G)$, which leads to the definition of the Frobenius digraph of the codegree graph of $G$.
\begin{defi}
Let $G$ be a finite group. The \emph{Frobenius digraph of codegree graph} of $G$, denoted $\overrightarrow{\Gamma}(G)$, is the directed graph obtained from $\overline{\Gamma}(G)$ by orienting each edge $pq$ as follows:
for the Hall $\{p,q\}$-subgroup $PQ \leq G$,
\begin{enumerate}
    \item if $PQ$ is a Frobenius group with kernel $Q$ and complement $P$, orient $pq$ as $p \to q$;
    \item if $PQ$ is a $(p,q,p)$-type 2-Frobenius group, also orient $pq$ as $p \to q$.
\end{enumerate}
\end{defi}

The orientation rule above follows the convention for the Frobenius digraph of the prime graph introduced in~\cite{gruber2015}. The resulting directed graph is denoted by $\overrightarrow{\Gamma}(G)$.
It remains to verify that this orientation is well-defined.

\begin{thm}
     The Frobenius digraph of the codegree graph is well-defined, thus it is independent of the choice of Hall $\{p,q\}$-subgroup.
\end{thm}

\begin{proof}
By Lemma \ref{lem:Q_leq_Sol(G)}, if $p \to q$ in $\overrightarrow{\Gamma}(G)$, then $Q \leq \Sol(G)$, where $\Sol(G)$ denotes the largest solvable normal subgroup of $G$. 
Suppose for contradiction that there exists another Hall $\{ p,q\}$-subgroup $P^\star Q^\star$ such that the orientation of edge $pq$ is  $q\to p$. Then $P^\star\leq \Sol(G)$. Since $\Sol(G) \trianglelefteq G$ and Sylow subgroups of $G$ are conjugate, it follows that $PQ \leq \Sol(G)$ and $P^\star Q^\star \leq \Sol(G)$. Moreover, since $\Sol(G)$ is solvable, all its Hall $\{p,q\}$-subgroups are conjugate, hence isomorphic. In particular  $PQ \cong P^\star Q^\star$. However, by construction of the orientation, isomorphic Frobenius or 2-Frobenius groups must induce the same orientation on $\{p,q\}$, this contradicts the assumption that $PQ$ induces $p\to q$ while $P^\star Q^\star$ induces $q\to p$.
Thus, the orientation of $\overrightarrow{\Gamma}(G)$ is well-defined and independent of the choice of Hall $\{p,q\}$-subgroup.
\end{proof}

\begin{lemma}\label{lem:codegree_no_3_digraph}
    The Frobenius digraph of a codegree graph $\overrightarrow{\Gamma}(G)$ contains no directed path of length $3$.
\end{lemma}

\begin{proof}

Suppose for contradiction that $a\to b \to c \to d$ is a directed 3-path in $\overrightarrow{\Gamma}(G)$. By Lemma \ref{lem:Q_leq_Sol(G)}, for any $B\in \Syl_b(G)$, $C \in \Syl_c(G)$, and $D \in \Syl_d(G)$, we have $BCD \leq \Sol(G)$. Thus $|\Sol(G)|_{\{b,c,d\}} = |G|_{\{b,c,d\}}$, where $|X|_\pi$ denotes the $\pi$-part of the order of $X$.

Let $\overline{X}\in \Syl_a(G/ \Sol(G))$, and let $X$ be its inverse image under the natural projection $G \to G/ \Sol(G)$, so  $X=\overline{X} \Sol(G)$. Since $\overline{X}$ and $\Sol(G)$ are both solvable, so is $X$. Therefore, $X$ contains a Hall $\{ a,b,c,d\}$-subgroup $H$ of $G$.  
Moreover, since $|X|_{\{ a,b,c,d\}}=|G|_{\{ a,b,c,d\}}$, $H$ is also a Hall $\{a,b,c,d\}$-subgroup of $G$. Since $H\leq X$ and $X$ is solvable, $H$ is also solvable.

Since the orientation of $\overrightarrow{\Gamma}(G)$ is independent of the choice of Hall subgroup, it follows that the directed path $a\to b \to c \to d$ is contained in $\overrightarrow{\Gamma}(H)$. It remains to show that $\overrightarrow{\Gamma}(H)$ contains no directed path of length $3$.
Let $\overrightarrow{\Gamma}_e(H)$ denote the Frobenius digraph of $H$ as defined in \cite{gruber2015}. By \cite[Corollary to Theorem E]{qian2007} the prime graph $\Gamma_e(H)$ is a subgraph of the codegree graph $\Gamma(H)$, which implies  $\overline{\Gamma}(H) \subseteq \overline{\Gamma}_e(H)$. 
Moreover, the orientation of each edge in $\overline{\Gamma}(H)$ is identical in both $\overrightarrow{\Gamma}(H)$ and $\overrightarrow{\Gamma}_e(H)$. Therefore, the directed $a\to b \to c \to d$ is contained in $\overrightarrow{\Gamma}_e(H)$.
But by \cite[Corollary 2.7]{gruber2015}, $\overrightarrow{\Gamma}_e(H)$ contains no directed $3$-path. A contradiction.
\end{proof}

\begin{proof}[Proof of Theorem 1.1]
The sufficiency follows from \cite[Proposition 3.5]{liu2025}, which establishes that if the complement $\overline{\Gamma}$ of a graph $\Gamma$ is triangle-free and $3$-colorable, then $\Gamma$ is realizable as the codegree graph of some finite solvable group $G$. 
We now prove the necessity. Let $\Gamma:=\Gamma(G)$ denote the codegree graph of a finite group $G$.
By Lemma \ref{lem:codegree_no_3_digraph} $\overrightarrow{\Gamma}(G)$ contains no directed $3$-path, from the Gallai-Roy Theorem  \cite[Theorem 7.17]{chromatic_graph_theory}, $\overline{\Gamma}(G)$ is $3$-colorable. Moreover by \cite[Theorem E]{qian2007}, $\overline{\Gamma}(G)$ is triangle-free. This completes the proof of necessity, and hence the theorem. 
\end{proof}

\section{Groups Whose Codegree Graph is a 5-Cycle}
In this section, we characterize the group whose codegree graphs is a $5$-cycle. As the application of the characterization of codegree graph, we borrow the conception of minimal prime graph  from \cite{gruber2015} to establish the concept of codegree graph. 
\begin{defi}
Let $\Gamma(G)$ be the codegree graph of group $G$. If $\Gamma(G)$ satisfies:
\begin{itemize}
    \item $|V(G)|>1,$
    \item $\Gamma(G)$ is connected,
    \item $\Gamma(G)-\{p,q\}$ is not the codegree graph for any $p,q\in V(G), $
\end{itemize}
then we call $\Gamma(G)$ is a \emph{minimal codegree graph}.
\end{defi}
Due to the identical graph-theoretic characterizations of the codegree graph and the prime graph, the minimal graphs share the same graph-theoretic properties. For example, $5$-cycle is the minimal codegree graph with the smallest number of vertices. And by \cite[Lemma 4.1]{gruber2015}, every minimal codegree graph contains an induced $5$-cycle.
 
To characterize the group whose codegree graph is $5$-cycle, we synthesize the findings presented in \cite[Proposition 3.5]{gruber2015}.
\begin{lemma}
     Let $G$ be a solvable finite group. Suppose that $\pi(G)=\{a,b,c,d,e\}$. Let $A,B,C,D,E$ be Sylow $a$-, $b$-, $c$-, $d$-, $e$-subgroups of $G$, respectively, such that $A,B,C$, $D,E$ form a Sylow system of $G$. Then, up to a permutation of the set $\pi(G)$, the prime graph $\Gamma_e(G)$ is a $5$-cycle if and only if one of the following conditions holds. 
     
\begin{itemize}
        \item [(1)] $G=\Fit(G) \rtimes (B\rtimes 
        AD)$ where $\Fit(G)=C \times E$, $ABC=\Frob_2(A,B,C)$, $DCE=\Frob(D,CE)$, $AE= \Frob(A,E)$.
        \item [(2)] $G=\Fit(G) \rtimes (EB\rtimes \overline{A}D)$ where $Fit(G)=C\times O_a(G)$, $\overline{A}\cong A/O_a(G)>1$, $ABC=\Frob_2 (A,B,C)$, $DCE=\Frob(D,CE)$, $AEB=\Frob_2(\overline{A},BE,O_a(G))$.
   \end{itemize}
\end{lemma}

\begin{proof}
    Up to isomorphism, there is a unique way to orient its Frobenius digraph $\overrightarrow{\Gamma}_e(G)$.
    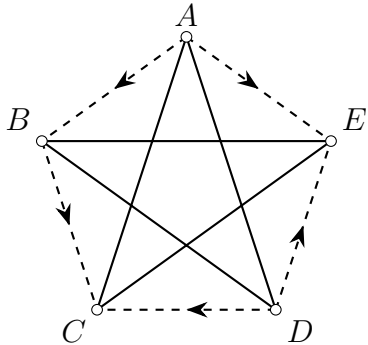
\begin{figure}[ht]
    \centering
    \begin{tikzpicture}[
    midarrow/.style={
        decoration={
            markings,
            mark=at position 0.5 with {\arrow{Stealth[scale=1.2]}}
        },
        postaction={decorate}
    }
]
        \coordinate (A) at (90:2cm);
        \coordinate (B) at (162:2cm);
        \coordinate (C) at (234:2cm);
        \coordinate (D) at (306:2cm);
        \coordinate (E) at (18:2cm);
       
       \draw [thick] (A)--(C)--(E)--(B)--(D)--cycle;

       \draw[thick, dashed, midarrow] (A)--(B);
       \draw[thick, dashed, midarrow] (A)--(E);
       \draw[thick, dashed, midarrow] (D)--(E);
       \draw[thick, dashed, midarrow] (D)--(C);
       \draw[thick, dashed, midarrow] (B)--(C);

        \filldraw[fill=white] (A) circle (2pt);
        \filldraw[fill=white] (D) circle (2pt);
        \filldraw[fill=white] (B) circle (2pt);
        \filldraw[fill=white] (E) circle (2pt);
        \filldraw[fill=white] (C) circle (2pt);

        \node[anchor=south] at (A) {$A$};
        \node[anchor=south east] at (B) {$B$};
        \node[anchor=north east] at (C) {$C$};
        \node[anchor=north west] at (D) {$D$};
        \node[anchor=south west] at (E) {$E$};
    \end{tikzpicture}
    \caption{5-cycle and its Frobenius digraph.}
    \label{fig:5-ccyle}
\end{figure}

    By \cite[Property 3.5]{gruber2015}, we have $C\leq \Fit(G)$. Since $b$ and $d$ are non-adjacent to $c$ in $\overline{\Gamma}_e(G)$, it follows $b,d \notin \pi(\Fit(G))$. Moreover since $a$ and $e$ are non-adjacent, at most one of $a$ or $e$ lies in $\pi(\Fit(G))$.

    Suppose $a\in (\Fit(G))$, $\Fit(G)=O_a(G)\times C$. By the Frobenius direction from $a$ to $b$ and $c$, it follows that  $AB$ and $AE$ are both 2-Frobenius groups with $\overline{A}\cong  A/ O_a(G)>1$, and by \cite[Property 4.4]{gruber2015}, $l_F(G)=3$, so $l_F(G/\Fit(G)) = 2$. Thus, the prime graph of $\overline{G}:=G/\Fit(G)$ is a path of length $3$ and the Fitting height of $G$ is $2$. By \cite[Theory 1.3]{qian2025diameter3},  $\overline{G} = EB\rtimes \overline{A}D$, so $G=\Fit(G) \rtimes(EB\rtimes \overline{A}D)$, where $\Fit(G)=C\times O_a(G)$. The Frobenius actions among these subgroups are uniquely determined by the orientation of the Frobenius digraph $\overrightarrow{\Gamma}(G)$.
    
    Suppose $a\notin \pi (\Fit(G))$. By \cite[Property 4.4]{gruber2015}, the Fitting length of $G$ is $l_F(G)=3$. it follows that $A\leq \overline{G}:=G/\Fit(G)$. Since the Fitting height of $\overline{G}$ is $2$, $AB$ is Frobenius group. By \cite[Property 3.4]{gruber2015}, $AE$ is Frobenius group and $BE=B\times E$. Similarly, since $d\notin \Fit(G)$ , $DE$ is a Frobenius group and by \cite[Property 3.4]{gruber2015}, $CE=C\times E$. Now for each prime $r\in \pi(G)$, there exists a Sylow $r$-subgroup normalizes $E$, and hence $E\trianglelefteq G$, in other words, $\Fit(G)=C\times E$. Then by \cite[Property 4.4]{gruber2015}, $l_F(G/\Fit(G))=2$. Thus, for $\overline{G}:=G/\Fit(G)$, it satisfies that $\overline{G}=B\rtimes AD$. Since $\Fit(G)=C\times E$, The Frobenius actions among these subgroups are uniquely determined by the orientation of Frobenius digraph $\overrightarrow{\Gamma}_e(G)$.
\end{proof}

 When $\Gamma(G)=C_5$, By \cite[Theorem 3.1]{alizadeh2021}, since $\Gamma(G)$ is triangle-free, $G$ is also a solvable group. 

\begin{proof}[Proof of Theorem 1.2]
First we prove necessity.
Suppose the codegree graph $\Gamma(G)$ of a group $G$ is a $5$-cycle. Moreover, by \cite[Theorem E]{qian2007}, the prime graph $\Gamma_e(G)$ is the subgraph of $\Gamma(G)$. However, by the property of the prime graph of a solvable group \cite[Theorem 2.10]{gruber2015}, no proper subgraph of $5$-cycle can be the prime graph of a solvable group, in other words, $5$-cycle is a minimal prime graph in the sense of \cite{gruber2015}. So $\Gamma_e(G)$ is also a $5$-cycle. Consequently, the complement graph of the prime graph, $\overline{\Gamma}(G)$ is also a $5$-cycle. 

We now prove sufficiency. 
Suppose $G$ satisfies condition (1), $\Fit(G)=C\times E$. For every $\chi \in \Irr (G)$, we analyze the restriction of  $\Fit(G)=C\times E$. 

If $C\nsubseteq ker(\chi)$ and $E \nsubseteq ker(\chi)$, since $B,D$ act fixed-point-freely on $C$, by Lemma \ref{lem:normal_sylow634}, $cod(\chi)$ is a $\{b,d\}'$-number.
If $E \nsubseteq ker(\chi)$, since $A$ and $D$ act fixed-point-freely on $E$, by Lemma \ref{lem:normal_sylow634}, we have $cod(\chi )$ is an $\{a,d\}'$-number. Therefore, $cod(\chi)$  is a $\{ c,e\}$-number.
Otherwise, if $E\subseteq ker(\chi)$, $cod(\chi)$ is a $e'$-number. So $cod(\chi)$ is an $\{a,c\}$-number.

If $C\subseteq ker(\chi)$, then $cod(\chi)$ is a $c'$-number.
If $E \nsubseteq ker(\chi)$, similarly, since $A$ and $D$ act fixed-point-freely on $E$, by Lemma \ref{lem:normal_sylow634}, we have $cod(\chi )$ is an $\{a,d\}'$-number. Then $cod(\chi)$ is a $\{b,e\}$-number.
Otherwise, if $E\subseteq ker(\chi)$, $cod(\chi)$ is a $e'$-number. and $\chi$ may be view as an irreducible character on $\overline{G}:=G/\Fit(G) \cong B\rtimes AD$. We now analyze $\lambda_B\in Irr(\chi_B)$. If $B\nsubseteq ker \chi$, analyze $I_{\lambda_B}(\overline{G})$, since $A$ act fixed-point-freely $B$, by Lemma \ref{lem:normal_sylow634}, $cod(\chi)$ is $\{b,d\}$-number. If $B\subseteq ker \chi$, $\chi$ may be view as the character on $AD$, and $cod(\chi)$ is an $\{a,d\}$-number.

In all cases, the set of prime divisors of $\cod(\chi)$ 
is contained in one of $\{c,e\}$, $\{a,c\}$, $\{b,e\}$, $\{b,d\}$, or $\{a,d\}$. Therefore, $E(\Gamma(G))\subseteq\{ce,ac,be,bd,ad\}$, by \Cref{thm:codegree_graph}, any proper subgraph can't be a codegree graph of a group. Hence, $\Gamma(G)$ is a $5$-cycle.

Suppose $G$ satisfies condition(2). Then $\Fit(G)=C\times O_a(G)$. As in the proof for condition(1), it is suffices to show that  $E(\Gamma(G))\subseteq\{ce,ac,be,bd,ad\}$. 

If $A\nsubseteq ker \chi$, we have $\cod(\chi)$ is a $\{b,e\}'$-umber. If $C\nsubseteq ker\chi$, then $\cod(\chi)$ is $\{b,d\}'$-number. And $\cod(\chi)$ is an $\{a,c\}$-number. If $C\subseteq ker\chi$, $\cod(\chi)$ is a $c'$-number, so $\cod(\chi)$ is an $\{a,d\}$-number.

If $O_a(G)\subseteq ker\chi$, $\chi$ may be view as a character on $\tilde{G}:=G/O_a(G)$. Since $\overline{A}>1$, the codegree graph $\Gamma(\tilde{G})$ remains a $5$-cycle. In $\tilde{G}$, from the Hall $\{a,b,e\}$-subgroup $\overline{A}BE=\Frob(\overline{A},BE)$, we know that $\overline{A}\leq N_{\tilde{G}}(E)$, and since the Frobenius kernel $BE$ is a nilpotent group, $BE=B\times E$, $B\leq N_{\tilde{G}}(E)$. Similarly, from Hall $\{d,c,e\}$- subgroup $DCE=\Frob(D,CE)$, we have $CD\leq N_{\tilde{G}}(E)$. So for every prime in $\pi (G)$, there exist some Sylow subgroup that normalizes $E$. $E\trianglelefteq G$. Thus, $\Fit (\tilde{G})=C\times E$. And in $\tilde{G}$, $\overline{A}BC=\Frob_2(\overline{A},BC)$, $DCE=\Frob(D,CE)$, $\overline{A}E=\Frob(A,E)$, $\tilde{G}=\Fit(\tilde{G}) \rtimes (B\rtimes AD)$. Therefore, $\tilde{G}$ satisfies condition(1). Therefore, $E(\Gamma(G))\subseteq\{ce,ac,be,bd,ad\}$. It follows that $\Gamma(G)$ is also a $5$-cycle in case(2). 
\end{proof}

\section{Further Comparisons Between the Codegree Graph and the Prime Graph }
While codegree graphs and prime graphs share several common properties, there also exist groups whose codegree graph is distinct from the prime graph.
Even for the non-trivial situation with the least number of vertices and edges,  $3$ vertices and $2$ edges, Qian created a counterexample of a solvable group whose codegree graph is different from prime graph in \cite[Example 2.2]{qian2025}. The following conclusion comes from it. After summarizing the steps the counterexample we can get some new propositions such that we can create some counterexamples more practically.
 \begin{proposition}
     If $G$ has a Frobenius Hall $\{p,q\}$-subgroup, then $p$ and $q$ are not incident in prime graph $\Gamma_e(G)$. Let $N$ is a normal subgroup of a finite group $G$, such that $(|G:N|,|N|)=1$. If there exists elements $a$, $b$ and irreducible character $\lambda \in \Irr(N)$ such that $a$, $b$ fix $\lambda$ while $|a|=p$, $|b|=q$. Then $\lambda$ determines an irreducible character $\chi \in G$ satisfying $pq| \frac{|G|}{\chi(1)}$. Furthermore , if $a$ and $b \notin \ker \chi$, then there exists $\chi\in \Irr(G)$, $pq\mid \cod(\chi)$. Specially, $p$ and $q$ are incident in codegree graph $\Gamma(G)$.
 \end{proposition}
 \begin{proof}
     By \cite[Theorem 8.16]{isaac_character}, since $(|G:N|,|N|)=1$, then $\lambda$ can uniquely expand to $I_G(\lambda)$. Assume that $\psi$ is the irreducible character be expanded by $\lambda$. By \cite[Theorem 6.11]{isaac_character}, there exists an irreducible character $\chi \in Irr(G)$ satisfying $\psi^G=\chi$. Now $$\frac{|G|}{\chi (1)}=\frac{|G|}{\psi (1)|G:I_G(\lambda)|}=\frac{|I_G(\lambda)|}{\psi (1)}.$$
     While by the definition of inertia subgroup $I_G(\lambda)$, element $a$ and $b \in I_G(\lambda)$, so $pq\mid |I_G(\lambda)|$. Assume that $a$ and $b$ are also $\notin ker(\chi)$, then $pq \mid \cod(\chi)$.
 \end{proof}
\begin{lemma}\label{lem:second_step_of_counterexample}
    Let V be the unique normal subgroup of $G$. If $G$ satisfies $(|G:V|,|V|=1$, and a Frobenius subgroup $PQ=\Frob(P,Q)$, where exists two nontrivial subgroups $Q_1, Q_2$, satisfies $C_V(Q_1)=1$ and $C_V(Q_2)>1$. Then there exists $p q\mid \cod(\chi)$.
\end{lemma}
\begin{proof}
    Let $V_1=C_V(Q_1)$, since $(|G:V|,|V|)=1$, we analyze the Frobenius group $PQ_1$ acting on $V_1$. By \cite[Theorem 15.16]{isaac_character}, $C_{V_1}(P)>1$. Since the non-identity is fixed with $PQ_1$. By \cite[Theorem 13.24]{isaac_character}, the action of $Q_1\rtimes P$ on $V$ and to $\Irr(V)$ is permutation isomorphism. Since $(|G:V|,|V|)=1$, by \cite[Theorem 8.16]{isaac_character}, $\lambda$ can be expended to an irreducible character $I_\lambda (G)$. Then since $Q\rtimes P \leq I_\lambda(G)$, by \cite[Theorem 6.11]{isaac_character}, there exists $\chi\in Irr(\lambda^G)$, such that $pq\mid \frac{|G|}{\chi(1)}$. And $V$ is the unique normal subgroup fo $G$, while $\chi_V\neq 1$, it shows that $\chi$ is faithful. Hence, $pq|\cod(\chi)$.
\end{proof}

 \begin{lemma}\label{lem:third_step_of_counterexample}
     Let finite group $G=V\rtimes (Q\rtimes P)$, such that $(|G:V|,|V|)=1$ and $PQ=\Frob (P,Q)$. The action of $PQ$ on $V$ is faithful and irreducible. In $G$, $Q$ is abelian but not cyclic and $V$ is a nontrivial elementary abelian group. And there exists $p_0\in P$ normalize $C_Q(V)$. Then there exists $p q\mid \cod(\chi)$.
 \end{lemma}
 \begin{proof}
     First, since the action of $PQ$ on $V$ is faithful, and $V$ is an elementary abelian group, then $V$ is the unique minimal normal subgroup of $G$.
     We analyze the irreducible action of $Q$ on $V$. Since $Q$ is an abelian group, thus the action of $Q/C_Q(V)$ on $V$ is a faithful irreducible from an abelian group to $V$. Then $Q/C_Q(V)$ is a cyclic group. However, $Q$ is not a cyclic group, thus $C_Q(V)>1$. Assign $Q:=C_Q(V)$.
     Then we analyze $N_G(Q_0)$. 
     Since $Q\leq N_G(Q_0)$, $V\cap N_G(Q_0)=C_V(Q_0)$, and there exists $p\in P$ such that $\langle p_0\rangle \leq N_G(Q_0)$. 
     Thus $N_G(Q_0)\geq C_V(Q_0) \rtimes (Q\rtimes\langle p_0\rangle)$. 
     Since $(|P|,|Q|)=1$, the action of $\langle p_0\rangle$ on $Q$ is completely irreducible. 
     While $Q_0$ is $P_0$-invariant, so $C_{V_1}(B_0)>1$, and $C_{V_0}(B_1)=1$. 
     Thus, $\chi \in \Irr(G)$, by Lemma \ref{lem:second_step_of_counterexample}, $pq\mid \cod(\chi)$.
 \end{proof}
\begin{proposition}
    There exists infinity groups whose prime graph is a $2$-path while its codegree graph is a triangle.
\end{proposition}
\begin{proof}
    Let $PQ=\Frob(P,Q)$ where $P\cong \langle a\rangle \cong \mathbb{Z}_4$ and $Q\cong \mathbb{Z}_q \times \mathbb{Z}_q$. Let $V$ be an elementary abelian $r$-group such that $PQ$ act faithfully and irreducible on $V$. Set $G=V \rtimes H$, then $G$ satisfies our proposition, where $p=2,q,r$ are three different primes.
    Since $H$ has a faithful and absolutely irreducible representation of dimension $4$, as $(|PQ|,|V|)=1$, then $PQ$ has an irreducible module over the field $\mathbf{F}_r$.
    Since $PQ$ is a Frobenius group and $G$ is not a Frobenius group or a $2$-Frobenius group, then the prime graph of $G$ is a $2$-path.
    Since the action of $PQ$ on $V$ is irreducible, then $C_Q(V)\neq Q$. And $\langle a^2 \rangle$ normalizes the proper subgroup of $Q$, thus it normalizes $C_Q(V)$. Finally by Lemma \ref{lem:third_step_of_counterexample}, there exists $\chi \in \Irr  (G)$, $pq \mid \cod(\chi)$. Since $\Gamma_e(G)$ is a subgraph of $\Gamma(G)$, then $\Gamma(G)$ is a triangle.
\end{proof}

The following example describes a nontrivial construction of a Frobenius digraph whose prime graph is equal to its codegree graph.
\begin{proposition}
    Let $G$ be a solvable group such that the Frobenius digraph of its prime graph is as shown in \Cref{fig:6-points example}. Then the codegree graph of $G$ coincides with its prime graph.
\end{proposition}
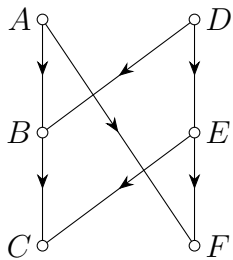
\begin{figure}[ht]
    \centering
    \begin{tikzpicture}[
    midarrow_2/.style={
        postaction={decorate},
        decoration={
            markings,
            mark=at position 0.5 with {\arrow{Stealth[scale=1.2]}}
        }
    }
]
        \coordinate (A) at (0cm,3cm);
        \coordinate (B) at (0cm,1.5cm);
        \coordinate (C) at (0cm,0cm);
        \coordinate (D) at (2cm,3cm);
        \coordinate (E) at (2cm,1.5cm);
        \coordinate (F) at (2cm,0cm);
       
       \draw [midarrow_2] (A)--(B);
       \draw [midarrow_2] (D)--(B);
       \draw [midarrow_2] (D)--(E);
       \draw [midarrow_2] (B)--(C);
       \draw [midarrow_2] (A)--(F);
       \draw [midarrow_2] (E)--(F);
       \draw [midarrow_2] (E)--(C);

        \filldraw[fill=white] (A) circle (2pt);
        \filldraw[fill=white] (D) circle (2pt);
        \filldraw[fill=white] (B) circle (2pt);
        \filldraw[fill=white] (E) circle (2pt);
        \filldraw[fill=white] (C) circle (2pt);
        \filldraw[fill=white] (F) circle (2pt);

        \node[anchor=east] at (A) {$A$};
        \node[anchor=east] at (B) {$B$};
        \node[anchor=east] at (C) {$C$};
        \node[anchor=west] at (D) {$D$};
        \node[anchor=west] at (E) {$E$};
        \node[anchor=west] at (F) {$F$};
    \end{tikzpicture}
    \caption{An example of Frobenius digraph whose prime graph coincides codegree graph}
    \label{fig:6-points example}
\end{figure}

\begin{proof}
    Since $diam(\overline{\Gamma}_e(G))=2$. By the definition in \cite{gruber2015}, $\overline{\Gamma}(G)$ is a minimal prime graph. It follows from \cite[Proposition 3.5]{gruber2015} that the Sylow $c$-subgroup $C$ and the Sylow $f$-subgroup $F$ are normal in group $G$. Now, consider any irreducible character $\chi \in \Irr (G)$, we analyze $\lambda_C \in \Irr (\chi_C)$. If $C\nleq ker\chi$, then since $B$ and $E$ act fixed-point-freely on $C$, by Lemma \ref{lem:normal_sylow634},
    $\cod(\chi)$ is a $\{b,e\}'$-number. Furthermore, since $F\trianglelefteq G$, and $A$ act fixed-point-freely on $F$, it follows from Lemma \ref{lem:normal_sylow634} that
    the edges incident to the vertex $f$ are identical in both the prime graph and the codegree graph. Hence, $\cod(\chi)$ is either an $a'$-number or an $f'$-number.
    Consequently, when $C\nleq ker \chi$, the codegree $\cod(\chi)$ is a $\{c,d,f\}$-number or an $\{a,c,d\}$-number. Assign  $\Gamma(G|C)$ as the Gruenberg-Kegel graph of codegrees of characters of $G$ satisfying their kernel doesn't include $C$. Then $\Gamma(G|C)\subseteq \Gamma_e(G)$.
    If $C\leq ker\chi$, then $\chi$ may be viewed as a character on $G/C$. Since $\Gamma_e$ is a $5$-cycle, we have $\Gamma(G/C)=\Gamma_e(G/C)$. As every edges in $\Gamma(G)$ corresponds to an edge in $\Gamma_e(G)$, it follows that $\Gamma(G)=\Gamma_e(G)$.
\end{proof}
It is readily observed that we can construct an arbitrary number of vertices that shares the same adjacency relationship with vertex $c$. Consequently, we can construct directed graphs of arbitrarily large size whose codegree graph coincides with their prime graph.
\begin{corollary}
    Let solvable group $G$ satisfy its Frobenius digraph of codegree like \Cref{fig:infinite point example}. 
    Then the codegree graph is same to its prime graph.
\end{corollary}
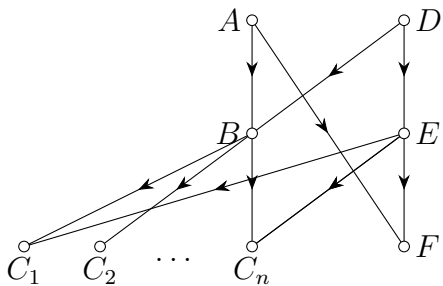
\begin{figure}[ht]
    \centering
    \begin{tikzpicture}[
    midarrow_2/.style={
        postaction={decorate},
        decoration={
            markings,
            mark=at position 0.5 with {\arrow{Stealth[scale=1.2]}}
        }
    }
]

        \coordinate (A) at (0cm,3cm);
        \coordinate (B) at (0cm,1.5cm);
        \coordinate (C) at (0cm,0cm);
        \coordinate (D) at (2cm,3cm);
        \coordinate (E) at (2cm,1.5cm);
        \coordinate (F) at (2cm,0cm);
        \coordinate (C_i) at (-1cm,0cm);
        \coordinate (C_2) at (-2cm,0cm);
        \coordinate (C_1) at (-3cm,0cm);

       \draw [midarrow_2] (A)--(B);
       \draw [midarrow_2] (D)--(B);
       \draw [midarrow_2] (D)--(E);
       \draw [midarrow_2] (B)--(C);
       \draw [midarrow_2] (A)--(F);
       \draw [midarrow_2] (E)--(F);
       \draw [midarrow_2] (E)--(C);
       \draw [midarrow_2] (B)--(C_1);
       \draw [midarrow_2] (B)--(C_2);
       \draw [midarrow_2] (E)--(C_1);
       \draw [midarrow_2] (E)--(C);

        \filldraw[fill=white] (A) circle (2pt);
        \filldraw[fill=white] (D) circle (2pt);
        \filldraw[fill=white] (B) circle (2pt);
        \filldraw[fill=white] (E) circle (2pt);
        \filldraw[fill=white] (C) circle (2pt);
        \filldraw[fill=white] (C_1) circle (2pt);
        \filldraw[fill=white] (C_2) circle (2pt);
        \filldraw[fill=white] (F) circle (2pt);

        \node[anchor=east] at (A) {$A$};
        \node[anchor=east] at (B) {$B$};
        \node[anchor=north] at (C) {$C_n$};
        \node[anchor=north] at (C_1) {$C_1$};
        \node[anchor=north] at (C_2) {$C_2$};
        \node[anchor=west] at (D) {$D$};
        \node[anchor=west] at (E) {$E$};
        \node[anchor=west] at (F) {$F$};
        \node[anchor=north] at (C_i) {$\dots$};
    \end{tikzpicture}
    \caption{An arbitrary vertices example of Frobenius digraph whose prime graph coincides codegree graph}
    \label{fig:infinite point example}
\end{figure}

The following lemma extends \cite[Proposition 3.5]{gruber2015} from the prime graph of a solvable group to the codegree graph of an arbitrary finite group.
\begin{proposition}\label{prop:minimal_2_dipath}
    Let $\Gamma(G)$ be a minimal codegree graph of finite group $G$. If a prime $p$ is an endpoint of a directed path length $2$ in the codegree Frobenius digraph $\overrightarrow{\Gamma}(G)$, then its corresponding Sylow $p$-subgroup $P$ is normal in $G$.
\end{proposition}
The proof of Proposition \ref{prop:minimal_2_dipath} is almost the same as in  \cite[Proposition 3.5]{gruber2015}. The only difference is that for every Hall $\pi$-subgroup, we invoke \ref{lem:Q_leq_Sol(G)} to prove it exists and that it lies in $\Sol(G)$.

The following result extends \cite[Theorem C(2)]{qian2025} from $(p,q,p)$-type 2-Frobenius group to a directed path of length $2$.
\begin{thm}
    Let $r,p,q\in\pi(G)$ and $R,P,Q$ are respectively Sylow $r$-, $p$-, $q$-subgroup in $G$. If there exists a directed path of length $2$ in $\overrightarrow{\Gamma(G)}$, which is $r \to p \to q$. If $|P|>p$, then $Q \trianglelefteq G$.
\end{thm}
\begin{proof}
    First, since $r \to p \to q$ in $\overrightarrow{G}$, by \ref{lem:Q_leq_Sol(G)}, $PQ\leq Sol(G)$. Then $G$ and $P^G$ are $\{p,q\}$-solvable. By \cite[Lemma 3.1]{qian2025}, $Q\leq P^G$ and $P^G=(PQ)^G$. Assign $D:=P^G=(PQ)^G$. Since $G$ is $p$-solvable and $P$ is cyclic, $l_p(G)=1$. Hence, $D$ has a normal $p$-complement. By \cite[Theorem C(1)]{qian2025}, $Q^D$ has a normal $p$-complement, $Q^G=Q^D\leq O_{p'}(D)$. Since $r\to p$, there exist $r$-group act fixed-point-freely on $P$. Let $G_1:=Q^G\rtimes(P\rtimes R_1)<G$, by \cite[Theorem B]{qian2025}, $r\to p \to q$ in $\overrightarrow{G_1}$, by induction, $Q^{G_1}=Q$. Since $Q^G<G_1$, then $Q\trianglelefteq Q^G$, and $Q$ is the only Sylow $q$-subgroup in $Q^G$, $Q \  \mathrm{char}\ Q^G$. And $Q\trianglelefteq G$, the conclusion is done in this case.
    Assume that $G=Q^G\rtimes(P\rtimes R_1)$. Then $D=P^G=(PQ)^G=Q^G\rtimes P$. Moreover, $Q^G=Q^D$ has a normal $q$-complement. Let $V\rtimes Q=Q^G$, where $V=O_{q'}(Q^G)=O_{\pi '}(Q^G)=O_{\pi '}(G)$.
    Since $P^G$ is solvable, and $R$ is cyclic, $G$ is solvable. To show $Q\trianglelefteq G$, we just need to show $V=1$. Let G be the minimal counterexample for $Q\trianglelefteq G$.

    Case 1: There exists minimal normal subgroup $E$ different from $V$. Then either $E<V=O_{\pi'}(G)$ or $E<Q$. Therefore, $\overline{G}:=G/E$ satisfies $r\to p \to q$ included in $\overrightarrow{\Gamma}(G)$. By induction, $\overline{Q}\trianglelefteq \overline{G}$. If $E<Q$, then $Q\trianglelefteq G$, a contradiction. If $E<V$, then $QE \trianglelefteq G$, then $Q^G \leq QE < QV=Q^G$, a contradiction.

    Case 2: The unique minimal normal subgroup of $G$ is $V$.
    Since $G$ is solvable, the order of $V$ is a prime power. Assume that $\pi(V)=r$. Then $P$ act fixed-point-freely on $V$ and $Q$, hence, $PVQ=\Frob(P,VQ)$. Now $VQ$ is a Frobenius kernel, so is nilpotent, then $V$ and $Q$ centralize each other. Specially, since $V\subseteq N_G(Q)$, $Q\trianglelefteq G$, a contradiction. Therefore, $r\notin \pi(V)$.
    Let $Z_p \cong B<P$, consider the co-prime action from $\Frob (P,Q)$ on $V$. Since $C_V(Q) \trianglelefteq G$, $V$ is the unique minimal normal subgroup and $Q\ntrianglelefteq G$, then $C_V(Q)=1$. By \cite[Theorem 15.16]{isaac_character}, $V_0:=C_V(B)>C_V(P)>1$. Now consider $N_G(B)$. Since  $BQ=\Frob(B,Q)$, $C_{QV}(B)=C_V(B)=V_0$, thus, $$N_{QV}(B)=C_{QV}(B)=C_V(B)=V_0.$$
    By $B\trianglelefteq Q_1P_1$, $$N_G(B)=N_G(B)\cap[(V\rtimes Q)\rtimes(P\rtimes R)]=N_{QV}(B)\rtimes(P\rtimes R)= V_0\rtimes(P \rtimes R).$$
    Consider the action from $\Frob(Q_1,P)$ to $V_0:=C_V(B)$. By the definition of $V_0$, the action from $P$ to $V_0$ is non-trivial. Since the action from $R_1P$ to $V_0$ is co-prime, then it is absolutely irreducible. Hence, there exists $V_1 \subseteq V_0$ such that $R_1P_1$ irreducibly act on $V_1$ and satisfies $C_{V_1}(P)=1$. By \cite[Theorem 15.16]{isaac_character}, $C_{V_1}(R_1)>1$, thus $B\rtimes R_1$ centralize a non-identity element in $V_1$. By \cite[Theorem 13.24]{isaac_character}, the action of $B\rtimes R_1$ to $V$ and to $Irr(V)$ is permutation isomorphic. Thus, there exist a non-principle $\lambda \in Irr(V)$ satisfying $\lambda$ is fixed by $B\rtimes R_1$. Since $(|G:V|,|V|)=1$, by \cite[Theorem 8.16]{isaac_character}, $\lambda$ can expand on $I_G(\lambda)$. By \cite[Theorem 6.11]{isaac_character}, there exists $\chi \in \Irr(\lambda^G)$ satisfying $pq\mid \frac{|G|}{\chi(1)}$. Since $V$ is the unique normal subgroup in $G$, then $\chi$ is faithful and $pq\mid \cod(\chi)$, a contradiction. 
\end{proof}

\vskip0.2in
\thanks{{\bf Acknowledgements.}
The authors are grateful to the referee for his/her valuable suggestions.}

\bibliographystyle{abbrv}
\bibliography{reference}	

\begin{thebibliography}{10}

\bibitem{alizadeh2020}
N.~Ahanjideh.
\newblock The {F}itting subgroup, {$p$}-length, derived length and character table.
\newblock {\em Math. Nachr.}, 294(2):214--223, 2021.

\bibitem{alizadeh2019}
F.~Alizadeh, H.~Behravesh, M.~Ghaffarzadeh, M.~Ghasemi, and S.~Hekmatara.
\newblock Groups with few codegrees of irreducible characters.
\newblock {\em Commun. Algebra}, 47(3):1147--1152, 2019.

\bibitem{alizadeh2021}
F.~Alizadeh, M.~Ghasemi, and M.~Ghaffarzadeh.
\newblock Finite groups whose codegrees are almost prime.
\newblock {\em Comm. Algebra}, 49(2):538--544, 2021.

\bibitem{chromatic_graph_theory}
G.~Chartrand and P.~Zhang.
\newblock {\em Chromatic graph theory}.
\newblock Discrete Mathematics and its Applications (Boca Raton). CRC Press, Boca Raton, FL, 2009.

\bibitem{du2016}
N.~Du and M.~L. Lewis.
\newblock Codegrees and nilpotence class of {{\(p\)}}-groups.
\newblock {\em J. Group Theory}, 19(4):561--567, 2016.

\bibitem{gruber2015}
A.~Gruber, T.~M. Keller, M.~L. Lewis, K.~Naughton, and B.~Strasser.
\newblock A characterization of the prime graphs of solvable groups.
\newblock {\em J. Algebra}, 442:397--422, 2015.

\bibitem{isaac_character}
I.~M. Isaacs.
\newblock {\em Character theory of finite groups}, volume No. 69 of {\em Pure and Applied Mathematics}.
\newblock Academic Press [Harcourt Brace Jovanovich, Publishers], New York-London, 1976.

\bibitem{liu2025}
Y.~Liu and Y.~Yang.
\newblock Two results on character codegrees.
\newblock {\em J. Algebra Appl.}, 24(6):Paper No. 2550158, 7, 2025.

\bibitem{qian2025}
G.~Qian.
\newblock Finite groups with non-complete character codegree graphs.
\newblock {\em J. Algebra}, 669:75--94, 2025.

\bibitem{qian2025diameter3}
G.~Qian.
\newblock Finite solvable groups whose prime graphs have diameter 3.
\newblock {\em Acta Math. Sin. (Engl. Ser.)}, 41(3):975--984, 2025.

\bibitem{qian2007}
G.~Qian, Y.~Wang, and H.~Wei.
\newblock Co-degrees of irreducible characters in finite groups.
\newblock {\em J. Algebra}, 312(2):946--955, 2007.

\bibitem{qian2002}
G.~H. Qian.
\newblock Notes on character degree quotients for finite groups.
\newblock {\em J. Math. (Wuhan)}, 22(2):217--220, 2002.

\bibitem{qian_overview}
G.~H. Qian.
\newblock Character codegrees in finite groups.
\newblock {\em Adv. Math. (China)}, 52(1):1--13, 2023.

\bibitem{yang2017}
Y.~Yang and G.~Qian.
\newblock The analog of {Huppert}'s conjecture on character codegrees.
\newblock {\em J. Algebra}, 478:215--219, 2017.

\end{thebibliography}
		
\end{document}